\documentclass[a4paper,11pt]{article}
\usepackage{array}
\usepackage{theorem}
\usepackage{amsmath,amscd,amssymb} 
\usepackage{latexsym}
\usepackage{epsfig}

\theoremstyle{plain}

\newtheorem{theorem}{Theorem}
\newtheorem{question}{Question}

\newcommand{\CC}{\mathbb C}

\newcommand{\NN}{\mathbb N}

\newcommand{\PP}{\mathbb P}

\newcommand{\RR}{\mathbb R}

\newcommand{\ZZ}{\mathbb Z}

\newcommand{\cD}{\mathcal D}

\newcommand{\cF}{\mathcal F}

\newcommand{\cH}{\mathcal H}

\newcommand{\tensor}{\otimes}
\newcommand{\emb}{\hookrightarrow}

\renewcommand{\Tilde}{\widetilde}

\renewcommand{\Bar}{\overline}

\newcommand{\cross}{\times}

\newcommand{\imic}{\cong}

\newcommand{\Orth}{\mathop{\null\mathrm {O}}\nolimits}

\newcommand{\Hom}{\mathop{\mathrm {Hom}}\nolimits}

\newcommand{\latt}[1]{{\langle{#1}\rangle}}

\newcommand{\Kthree}{\mathop{\mathrm {K3}}\nolimits}
\newcommand{\qedsymbol}{\mbox{$\Box$}}
\newcommand{\qed}{\unskip\nobreak\hfil\penalty50\hskip1em\hbox{}\nobreak
\hfill\qedsymbol\parfillskip=0pt\finalhyphendemerits=0}

\newenvironment{ProofwCaption}[1]
 {\addvspace\theorempreskipamount \noindent{\it #1.}\rm}
 {\qed \par \addvspace\theorempostskipamount}

\begin{document}
\title{On some lattice computations related to moduli problems}
\author{A.~Peterson and G.K.~Sankaran, with an appendix by V.~Gritsenko}
\maketitle
\begin{abstract}
\noindent The method used in \cite{K3} to prove that most moduli
spaces of $\Kthree$ surfaces are of general type leads to a
combinatorial problem about the possible number of roots orthogonal to
a vector of given length in $E_8$. A similar problem arises for $E_7$
in \cite{sympl}. Both cases were solved partly by computer methods. We
use an improved computation and find one further case, omitted from
\cite{K3}: the moduli space $\cF_{2d}$ of $\Kthree$ surfaces with
polarisation of degree $2d$ is also of general type for $d=52$. We
also apply this method to some related problems. In Appendix A,
V. Gritsenko shows how to arrive at the case $d=52$ and some others directly.
\end{abstract}

\noindent Many moduli spaces in algebraic geometry can be described as
locally symmetric varieties, i.e.\ quotients of a Hermitian symmetric
domain $\cD$ by an arithmetic group $\Gamma$. One method of
understanding the birational geometry of such quotients is to use
modular forms for $\Gamma$ to give information about differential
forms on $\Gamma\backslash \cD$. In \cite{K3} this method was used to
prove that the moduli space $\cF_{2d}$ of polarised $\Kthree$ surfaces
of degree~$2d$ is of general type in all but a few cases. The method
works if there exists a modular form of sufficiently low weight with
sufficiently large divisor. In \cite{K3}, and again in \cite{sympl}
where a similar method was applied to certain moduli of polarised
hyperk\"ahler manifolds, the required modular form is constructed by
quasi-pullback of the Borcherds form $\Phi_{12}$.

A suitable quasi-pullback exists if a combinatorial condition is
satisfied: there should exist a vector $l$ in the root lattice
$E_8$ (or $E_7$ in the hyperk\"ahler case) of square $2d$, orthogonal
to very few roots. This is evidently the case if $d$ is large, but for
small~$d$ the search for such an $l$ invites the use of a
computer. This was done in both \cite{K3} and in \cite{sympl} by a
randomised search, relying on the large Weyl group to ensure that in
practice no cases would be missed. 

Here we present an exhaustive search carried out by the first
author. For the hyperk\"ahler case the exhaustive search confirmed the
results of the earlier randomised search, but in the $\Kthree$ case one
previously overlooked value of $d$ with a suitable vector was found,
namely $d=52$. In fact it turned out that the randomised search had
indeed found this value, and the omission of the case $d=52$ from
\cite{K3} happened because the output had been interpreted
incorrectly.
{\footnote{By me. -- GKS}}

Nevertheless the following result is true and has not previously
appeared in the literature.

\begin{theorem}\label{thm52} The moduli space $\cF_{2\cdot 52}$ of
$\Kthree$ surfaces with polarisation of degree $104$ is of general
type.
\end{theorem}
The paper is organised as follows. In Section~\ref{background} we
explain briefly what the combinatorial problem is and how it arises,
and give some more general combinatorial problems of the same
nature. In Section~\ref{computations} we describe the theoretical and
computational methods used to solve it, along with some other results
obtained in the same way. In Appendix~A, Valery Gritsenko explains
explains how the case $d=52$ could have been foreseen without the help
of a computer. Some of the computer code is given in Appendix~B.

\noindent\emph{Acknowledgements:} Part of this paper forms part of
A. Peterson's Masters' thesis. He would like to thank Gerard van der
Geer for his supervision, and the University of Amsterdam for the nice
environment it provides. The second author would like to thank the
Fondazione Bruno Kessler in Trento and the Max-Planck-Institut f\"ur
Mathematik in Bonn for support, and Valery Gritsenko for helpful
conversations.

\section{Combinatorial problems and moduli}\label{background}

In this section we first give a list of combinatorial questions and
then explain the geometry that originally motivated them. First we fix
some terminology. We say that $L$ is a \emph{lattice of signature
  $(a,b)$} if $L\imic\ZZ^{a+b}$ and we fix a bilinear form $(\,
,\,)\colon L\cross L\to\ZZ$ of signature $(a,b)$. If $x\in L$ we refer
to $(x,x)$ as $x^2$ and call it the \emph{length} of $x$. If the
length of $x$ is $2$ then $x$ is called a \emph{root}. If the roots
of $L$ generate $L$ as an abelian group then $L$ is called a
\emph{root lattice}. A lattice $L$ is \emph{unimodular} if it is equal
to its dual $L^\vee=\Hom(L,\ZZ)\supseteq L$. We do not assume that $L$
is always unimodular but for simplicity we do assume that $L$ is
\emph{even}, i.e.\ that $x^2$ is always an even integer.

$E_8$ denotes the unique even unimodular positive-definite lattice of
rank $8$, i.e.\ with signature $(8,0)$: this is the sign convention
of~\cite{Bou} and is also used in~\cite{K3}. If
$n\in 2\ZZ$ then $\latt{n}$ is the rank~$1$ lattice spanned by a
vector of length $n$, and $U$ denotes the integral hyperbolic plane
$\ZZ e+\ZZ f$ with $e^2=f^2=0$ and $(e,f)=1$. The symbol $\oplus$
denotes the orthogonal direct sum of lattices. If $\Lambda$ is a
lattice and $n\in \ZZ$, then $\Lambda(n)$ denotes the same lattice
with the quadratic form multiplied by~$n$. In particular, $E_8(-1)$ is
the negative-definite even unimodular lattice of rank~$8$.

\subsection{Combinatorial problems}\label{questions}

Let $\Lambda$ be a root lattice (usually it will be $E_8$ or $E_7$)
and denote by $R(\Lambda)$ the set of its roots,
i.e.\ $R(\Lambda)=\{r\in\Lambda\mid r^2=2\}$. The combinatorial
questions arising in \cite{K3} and \cite{sympl} are special cases of
the following.

\begin{question}\label{Qpq} Given integers $p>q\ge 0$, what are
the values of $d$ for which every vector of length $2d$ that is
orthogonal to at least $2q$ roots is orthogonal to at least $2p$
roots?
\end{question}

More generally we may ask about all possibilities.

\begin{question}\label{QP} 
Given an even natural number $2d$, what are the possible numbers of
roots orthogonal to a vector of length $2d$?
\end{question}

If $l\in\Lambda$ we denote by $R(l^\perp)$ the system of roots
of $\Lambda$ orthogonal to~$l$. We denote the answer to
Question~\ref{QP} by $P(\Lambda,d)$: that is
\begin{equation}\label{Pd}
P(\Lambda,d):=\{m\in\ZZ\mid \exists l \in\Lambda\ l^2=2d,\ \#
R(l^\perp)=m\}. 
\end{equation}
Thus $P(\Lambda,d)$ is a finite set of even non-negative integers. We
call this the \emph{root type} of the non-negative even integer $2d$
for the lattice $\Lambda$

There are some immediate restrictions on what the root type can be:
for example, if $\Lambda=E_8$ then the largest $m$ that can occur is
$126$, when $R(l^\perp)\imic E_7$; but in that case $l\in
(E_7)^\perp_{E_8}\imic A_1$, so $d$ must be a square.

Especially for $\Lambda=E_8$, the value of $m_0(d)=\min P(E_8,d)$ is
of interest as it determines the lowest weight of modular form
obtained by quasi-pullback (see Equation~\eqref{qpb} below). If
$m_0(d)=0$ then this form will not be a cusp form, so the value of
$m_1(d)=\min P(E_8,d)\cap\NN$ is also significant. We should also like
to know whether this form is unique. So we also have the following
questions.
\begin{question}\label{Qk0}
For given $d$ and $\Lambda$, how can we compute $m_0(d)$?
\end{question}
\begin{question}\label{Qlowest}
For given $m$, what is the smallest value $d(m)$ of $d$ for which
$m_1(d)\le m$?
\end{question}
If in Question~\ref{Qlowest} we replace $m_1$ by $m_0$, then the case
$m=0$ asks for the length of shortest vectors in the interior of a
Weyl chamber: these are the Weyl vectors, which are well known.

If $m\in P(\Lambda,d)$ there is a further natural refinement.
\begin{question}\label{Qorbits}
How many Weyl group orbits of vectors $l$ with $l^2=2d$ and $\#
R(l^\perp)=m$ are there?
\end{question}
Some values of $m$ are of particular interest for geometric reasons:
for instance, if $14\in P(E_8,d)$ then quasi-pullback of $\Phi_{12}$
gives a canonical form on $\cF_{2d}$ (see Section~\ref{geometry}
below). This leads us to the following variant of Question~\ref{Qpq}.
\begin{question}\label{Qroots} 
For given $m$ and $\Lambda$, what are the values of $d$ such that
$m\in P(\Lambda,2d)$?
\end{question}
We can compute the answers to some cases of these questions by the
methods described in Section~\ref{computations}.

\subsection{Moduli}\label{geometry}

The following construction describes several moduli spaces in algebraic
geometry, including the moduli of polarised $\Kthree$ surfaces.

Let $L$ be an even lattice of signature $(2,n)$. The Hermitian
symmetric domain associated with $L$ is $\cD_L$, one of the two
connected components of
\begin{equation*}
\cD_L\cup\Bar{\cD}_L=\{[w]\in \PP(L\tensor\CC)\mid w^2=0,\ (w,\Bar w)>0\}.
\end{equation*} 
The group $\Orth(L)$ of isometries of $L$ acts on this union and we
denote by $\Orth^+(L)$ the index~$2$ subgroup preserving $\cD_L$. The
action is discontinuous, with finite stabilisers, so if $\Gamma$ is
any finite index subgroup of $\Orth^+(L)$ then
\begin{equation*}
\cF_L(\Gamma):=\Gamma\backslash\cD_L
\end{equation*} 
is a complex analytic space. In fact it is a quasi-projective variety,
having a minimal projective compactification, the Baily-Borel
compactification $\cF_L(\Gamma)^*$, obtained by adding finitely many
curves (called $1$-dimensional cusps) meeting at finitely many points
($0$-dimensional cusps). It is often preferable to work with a
toroidal compactification $\Bar{\cF_L(\Gamma)}$, which is a modification of
$\cF_L(\Gamma)^*$ depending on some combinatorial choices at the
$0$-dimensional cusps.

A modular form for $\Gamma$ of weight $k$ and character $\chi\colon
\Gamma\to \CC^*$ is a holomorphic function $F$ on the affine cone
$\cD_L^\bullet\subset L\tensor \CC$ such that
\begin{equation*}
F(tZ)=t^{-k}F(Z)\ \ \forall t\in\CC^*\quad\text{ and
}\quad F(gZ)=\chi(g)F(Z)\ \ \forall g\in\Gamma.
\end{equation*}
$F$ is a cusp form if it vanishes at every cusp. For the cases we shall
consider the only possible characters are $1$ and $\det(g)$, and the
order of vanishing at a cusp is an integer: see~\cite{comms}.

The aim of \cite{K3} is to show that the moduli space $\cF_{2d}$ of
polarised $\Kthree$ surfaces of degree $2d$ is of general type for
most values of $d\in\NN$. Using the Torelli theorem for $\Kthree$
surfaces one can show that
\begin{equation*}
\cF_{2d}=\cF_{L_{2d}}(\Tilde\Orth^+(L_{2d})),
\end{equation*}  
where $\Tilde\Orth^+(L)$ is the finite index subgroup of $\Orth^+(L)$
that acts trivially on the discriminant group $L^\vee/L$ and
\begin{equation*}
L_{2d}:=2U\oplus 2E_8(-1)\oplus \latt{-2d}.
\end{equation*}
Modular forms of suitable weight can be interpreted as differential
forms on the moduli space provided that they have sufficiently large
divisor. Therefore, to prove that the moduli space is of general type
it is enough to give a sufficient supply of such modular forms. There
are several technical difficulties here, one of which is the presence
of singularities. A sufficient condition, however, was given
in~\cite{K3}.

\begin{theorem}\label{gtthm}
Suppose that $n\ge 9$ and that there exists a nonzero cusp form $F_a$
of weight $a<n$ and character $\chi\equiv 1$ or $\chi(g)=\det(g)$,
vanishing along any divisor $\cH\subset \cD_L$ fixed by reflections in
$\Gamma$. Then $\cF_L(\Gamma)$ is of general type.
\end{theorem}

The form $F_a$ is then used to give many forms of high weight with
sufficiently large divisor, of the form $F=F_a^kF_{(n-a)k}$, and these
in turn give pluricanonical forms on a smooth model of
$\Bar{\cF_L(\Gamma)}$.

To apply this in specific cases such as $\cF_{2d}$ one must therefore
construct~$F_a$. The method used in \cite{K3} to do this is
quasi-pullback of the Borcherds form $\Phi_{12}$. This construction
first appeared in~\cite{BKPS}. The Borcherds form itself was
constructed in~\cite{Borth} by means of a product expansion, whereby
its divisor is evident. It is a modular form (not a cusp form) of
weight~$12$ and character $\det$ for the group
$\Orth^+(II_{2,26})$. The lattice $II_{2,26}$ of signature $(2,26)$ is
$2U\oplus N(-1)$, where $N$ is any one of the 24
Niemeier lattices, positive definite unimodular lattices of rank~$24$:
see~\cite{CS}. For our purposes the correct choice of $N$ is $3E_8$. A
choice of a (not necessarily primitive) vector $l\in E_8$ of length
$2d$ gives an embedding
\begin{equation*}
L_{2d}=2U\oplus 2E_8(-1)\oplus \latt{-2d}\emb II_{2,26}=2U\oplus 3E_8(-1)
\end{equation*}
which in turn gives an embedding
\begin{equation*}
\cD_{L_{2d}}^\bullet\emb \cD_{II_{2,26}}^\bullet.
\end{equation*}
Denote the images of these embeddings by $L_{2d}[l]$ and $\cD^\bullet[l]$ respectively.

If $r\in L$ is a root it determines a Heegner divisor
$\cH^\bullet_r\subset \cD^\bullet_L$, given by the equation
$(Z,r)=0$. The Borcherds form vanishes (to order~$1$) along all the
Heegner divisors for $L=II_{2,26}$ and in particular its restriction
to $\cD^\bullet[l]$ vanishes, as needed to apply
Theorem~\ref{gtthm}. However, $\Phi_{12}\vert_{\cD^\bullet[l]}$ may
well be zero, since if $r$ is a root of $II_{2,26}$ orthogonal to
$L_{2d}[l]$ then $\cD^\bullet[l]\subset\cH^\bullet_r$.

Instead we take the quasi-pullback, simply dividing by the equation of
each such $\cH^\bullet_r$, noting that
$\cH^\bullet_{-r}=\cH^\bullet_r$. We put
\begin{equation*}
R_l=\{r\in R(II_{2,26})\mid (r,L_{2d}[l])=0\}\imic
\{r\in R(E_8)\mid (r,l)=0\}
\end{equation*}
and define the quasi-pullback to be
\begin{equation}
\label{qpb}
F[l]=\left.{\frac{\Phi_{12}}{\prod_{\pm r\in
      R_l}(r,Z)}}\right\vert_{\cD^\bullet[l]}.
\end{equation}
This is a nonzero modular form, and one can show that it is a cusp
form provided $R_l\neq \emptyset$. It vanishes along all the
Heegner divisors fixed by reflections in $\Orth^+(L_{2d})$.

The weight, however, goes up by~$1$ every time we divide, so the
weight of $F[l]$ is $12+\frac{1}{2}\# R_l$. We can therefore
show that $\cF_{2d}$ is of general type if we can find an $l\in
E_8$ of length $2d$ with $2\le \# R_l<2(n-12)=14$. Moreover, if we
can find a cusp form of weight precisely $n=19$ then, by a result of
Freitag \cite{Fr}, $\cF_{2d}$ has $p_g>0$ and in particular is not
uniruled.

This leads us to Question~\ref{Qpq}, with $q=1$ and $p=7$ or $p=8$,
for $\Lambda=E_8$. In~\cite{sympl}, similar considerations about the
moduli of some hyperk\"ahler manifolds with a certain type of
polarisation lead to Question~\ref{Qpq} with $q=1$ and $p=6$ or $p=7$,
for $\Lambda=E_7$.

\section{Solving the combinatorial problems}\label{computations}

The specific combinatorial problems encountered in \cite{K3} and
\cite{sympl} can be solved in principle by first bounding $d$. It is
clear that for sufficiently large $d$ an $l$ will exist orthogonal
to a number of roots in the required range: indeed, for sufficiently
large $d$ we can find $l$ orthogonal to exactly two roots. An
explicit bound, followed by a finite calculation, will solve the
problem. Neither is entirely straightforward, though. In \cite{K3} a
counting argument is used to show that an $l\in E_8$ with
$l^2=2d$, orthogonal to at least two and at most $12$ roots, exists
(and therefore $\cF_{2d}$ is of general type) unless
\begin{equation}\label{K3cond}  
28N_{E_6}(2d)+63N_{D_6}(2d)\ge 4N_{E_7}(2d),
\end{equation}
where $N_L(2d)$ is the number of ways of representing $2d$ by the
quadratic form~$L$. The inequality~\eqref{K3cond} certainly fails
for large $d$, but to obtain an effective bound on $d$ one must bound
$N_{E_6}(2d)$ and $N_{D_6}(2d)$ from above and $N_{E_7}(2d)$ from
below by explicit functions. This is a non-trivial problem in analytic
number theory but it can be done, and after some refinements it gives
a reasonable bound of around $d=150$. It would be possible to resort
to direct computation at that point, but there is no need yet. Some
integers in that range are excluded from the list of possibly
non-general type polarisations because the inequality~\eqref{K3cond}
(or another similar inequality) in fact fails. Others can be excluded
by inspection, actually producing a vector $l$ by guessing the root
system $R(l^\perp_{E_8})$. The root systems used in this way in
\cite{K3} were $4A_1$, $2A_1\oplus A_2$, $A_3$ and $A_1\oplus
A_2$. The root systems $3A_1\oplus A_2$ and $2A_2$ were not tried: see
Appendix~A.

In \cite{sympl} a similar procedure was used, although there is an
extra difficulty caused by the opposite parity of the rank: working in
$E_7$, one needs to estimate $N_R(2d)$ from above for some odd-rank
root systems $R$, and this problem is not so well studied as in the
even rank case.

In either case, eventually one is left with a residual list of values of $d$
for which the problem has not been settled. In \cite{K3} it consists
of most integers between $15$ and $60$ (for very small $d$ the moduli
space is known to be unirational). The residual problem in the
hyperk\"ahler case considered in \cite{sympl} is much smaller.

Now, if we want to be (reasonably) sure that no cases have been
missed, we do need a computer. Moreover, the methods we now use to
solve this problem can also be used to give answers to question such
as those posed in Section~\ref{questions}.

\subsection{Algorithms}\label{algorithms}

We begin by representing $E_8$ in the
usual way, as the set of points $l=(l_1,\ldots,l_8)\in \RR^8$ such
that the $l_i$ are either all integers or all strict half-integers
(i.e.\ either $l_i\in\ZZ$ for all $i$ or $2l_i$ is an odd integer for
all $i$) and $\sum l_i\in 2\ZZ$, with the standard Euclidean quadratic
form on~$\RR^8$.

We need a very rough upper bound on $N_{E_8}(2d)$, 
because we want to know whether $N_{E_8}(2d)$ is small enough to allow
a brute-force search for $l\in E_8$ with $l^2=2d$ having $2\le \#
R(l^\perp)\le 12$. We can easily find such a bound by noting that if
$l^2=2d$ then each of the $8$ components $l_i$ of $l$ must have $l_i^2
\leq 2 d$, so $-\sqrt{2 d} \leq l_i \leq \sqrt{2 d}$, and must be a
half-integer: that gives
\begin{equation}
  N_{E_8}(2d) \leq (2 \lfloor{2 \sqrt{2 d}}\rfloor + 1)^8
\end{equation}
For $d = 52$, this bound is about $8 \cdot 10^{12}$.

If we are a bit more precise, and note that the components of $l$ are
either all integers, or all proper (i.e.\ non-integer) half-integers,
we save a factor $2^7$, giving a bound of about $5 \cdot 10^{10}$. This is
within reach of a brute-force search, but it is still high, especially
considering that we have to do some substantial work for each
candidate (compute the inner product with 240 different
vectors\footnote{We can be a lot more efficient than that, and skip
  most of these inner products, but even then we still have to compute
  dozens of inner products per candidate vector.}).

Thus an exhaustive search of all vectors in $E_8$ of length $\le 60$
is not computationally impossible but it would be cumbersome and would
not extend to even slightly larger problems such as other cases of
Question~\ref{Qpq}. The Weyl group $W(E_8)$ has order $2^{14}\cdot 3^5\cdot
5^2\cdot 7=696729600$ and should be used to reduce the size of the
problem. There are two approaches to doing this.

\bigskip
\noindent{\bf A. Randomised search.} This is what was actually done in
\cite{K3} and \cite{sympl}. Since the non-existence of a vector $l$
gives no information about the moduli space, we are willing to accept
a very small probability of failing to detect such a vector. We
therefore choose a large number of vectors of length less than $2\cdot
61$ at random and expect that, as the Weyl group orbits are large,
every orbit will be represented.

This approach worked very fast, using only a laptop computer and
immediately available software (Maple). A search of twenty thousand
randomly chosen vectors found all the pairs $(d,\#R(l^\perp))$ in the
ranges wanted within the first two thousand iterations, in
approximately two minutes. That is fairly convincing practical
evidence that there are no more. Unfortunately the output was then
mistranscribed, leading to the omission of the case $d=52$ and the
erroneous (but not really misleading) statement in \cite{K3} that ``an
extensive computer search for vectors orthogonal to at least $2$ and
at most $14$ roots for other $d$ has not found any''.

It is noteworthy that a similar search in the case $\Lambda=E_7$ did
find some cases not discovered analytically, and for which a
constructive method of finding $l$ is still not known. In other
words, some cases of the main theorem of \cite{sympl} still have only
a computer proof, although once $l$ has been found it is easy
enough to verify its properties by hand.

It is not so easy to estimate the probability \emph{a priori} that a
Weyl orbit might be missed. The Weyl group of $R(l^\perp)$, which is a
subgroup of the Weyl group of $E_8$, obviously stabilises $l$ and has
order no more than $24$ if $\# R(l^\perp)\le 12$, but in principle the
stabiliser of $l$ in $W(E_8)$ could be much larger. In that case the
Weyl group orbit would be small and more easily missed. In practice
the randomised method seems to find all the orbits.
%
\bigskip

\noindent{\bf B. Exhaustive search.} The first author organised an
exhaustive search, exploiting the Weyl group by searching a
fundamental domain for the subgroup $H<W(E_8)$ generated by
permutations of the eight components $l_i$ and sign changes of an even
number of components. This subgroup $H$ has size $2^7 \cdot 8!$, so
index $135$ in $W(E_8)$: it gives us most of the symmetries, with very
little effort.

We say that $l \in E_8$ is in \emph{normal form} if its
components are all nonnegative (except possibly the first, $l_1$) and
the squares of the components are nondecreasing from low index to high
index. By acting with an element of $H$, we can translate any $l \in
E_8$ to one in normal form: first permute the components, so their
squares are in order; then make them all (but $l_1$) nonnegative, by
changing the sign of every negative component (except $l_1$), and
flipping the sign of $l_1$ once for every such change.

It is straightforward to enumerate the elements of length $2d$ in
$E_8$ that are in normal form. For brevity, we will describe this
only for the ones having integer components (one can get the ones with
proper half-integer components in a very similar manner).

\noindent\emph{Step 1.} For every index $i \neq 1$, in descending order, we
consider all the possible values of $l_i$: we require $l_i$ to be a
non-negative integer  such that
      \begin{itemize}
        \item its square, added to the sum of the squares of the
          coordinates that have been chosen (i.e.\ the $l_j^2$ with $j
          > i$), does not exceed $2 d$ (otherwise $l^2 > 2 d$, for any
          further choice of coordinates); and
        \item (unless $i = 8$) it is not greater than $l_{i + 1}$
          (otherwise $l$ would not be in normal form).
      \end{itemize}
In other words, we let $l_i$ take any value $s \in \ZZ$ such that
\begin{equation}
0\le s \leq \min\left\{l_{i + 1},\sqrt{2 d - \sum\nolimits_{j > i}
  l_j^2}\right\}.
\end{equation}
\noindent\emph{Step 2.} See if $2 d - \sum_{j = 2}^{8} l_j^2$ is a
perfect square $m^2$. If so, let $l_1$ take values $-m$ and $m$; if
not, discard this choice of coordinates.

\noindent\emph{Step 3.} Check whether the $l$ so obtained are in
$E_8$, i.e.\ whether $\sum_{j=1}^8 l_j \in 2 \ZZ$. Discard any that
are not in $E_8$.

We must then filter these enumerated $l\in E_8$ to find the ones
with $\# R(l^\perp)$ in the required range ($2\le \# R(l^\perp)\le 12$
for the case considered in~\cite{K3}): this part of the procedure is
exactly the same as for the randomised version. Since the roots come
in pairs $\pm r$ it is enough to take inner products with a prepared
list of positive roots ($120$ or them), and of course we can stop examining $l$ as soon as we find a seventh pair of roots orthogonal to it.

The first author implemented this search in a high-level programming
language (Haskell). Without spending much time optimising, this runs
fast enough (a second or so on commercial hardware, for each of
the low values
of $d$ we are interested in, namely $d\le 60$). The partial use of the
symmetries of $E_8$ is crucial, though: to go through all the vectors of given
length $2d$ would have taken weeks or months for a single value of~$d$.

This program discovered the lost case $d=52$ and therefore
Theorem~\ref{thm52}.  A variant of it for $E_7$ reconfirmed the
results obtained by the randomised method in~\cite{sympl}. The code used for
the $E_8$ case is given in Appendix~B.

\subsection{Further results}\label{further}

The exhaustive algorithm~(B) from Section~\ref{algorithms} can be
modified to compute, in reasonable time, answers to some of the
questions from Section~\ref{questions} for small values of the
parameters. We investigated Question~\ref{QP} and 
and Question~\ref{Qroots} for small $m$ and $d$ with $\Lambda=E_7$ and
$\Lambda=E_8$. For $\Lambda=E_8$ we also investigated
Question~\ref{Qorbits} for the particular case $m=14$, corresponding
to canonical forms on $\cF_{2d}$.

Specifically, we have so far computed the root type $P(\Lambda,2d)$
for $\Lambda=E_7$ and $\Lambda=E_8$ and $d\le 150$, and the first part
of the root type (whether $m\in P(\Lambda,2d)$ for $2\le m \le 20$,
say) for larger $d$, up to about $300$ (further for some values
of~$d$). This part of the computation is fairly fast and only minor
changes to the program are needed.  

A little more work, and more computer time, is needed for
Question~\ref{Qorbits}. We must work now with $W(E_8)$, not with $H$,
and we first compute a transversal for $W(E_8):H$ (representatives for
each of the $135$ left cosets of $H$) and then reduce each of the
$135$ translates of each $l$ to standard form before comparing them.

The outcome counts the number of ways of obtaining a canonical form on
$\cF_{2d}$ by quasi-pullback of $\Phi_{12}$. There is no assurance
either that the forms so obtained are linearly independent or that
there are not more canonical forms that do not arise this way. The
results are nevertheless intriguingly unpredictable. There are no such
vectors for $d<40$. There is such a vector for $d=40$, and also for
$d=42,\ 43,\ 48$ (two orbits), $49,\ 51$--$54,\ 55$ and $56$ (two
orbits each), $57$ and $59$. There is no such vector for $d=60$, but
for $61$ there are three orbits and thereafter the number of orbits
drifts upwards irregularly. Without further comment, we tabulate
below the number $\nu_{14}$ of $W(E_8)$ orbits of length
$2d$ vectors in $E_8$ orthogonal to exactly~$14$ roots for $61\le
d\le 150$.
\begin{center}\label{orbitstable}
\begin{tabular}{|c|c||c|c||c|c||c|c||c|c||c|c|}
\hline
$d$&$\nu_{14}$&$d$&$\nu_{14}$&$d$&$\nu_{14}$&$d$&$\nu_{14}$&$d$&$\nu_{14}$&$d$&$\nu_{14}$\\
\hline
61&3&76&1&91&5&106&2&121&4&136&8\\
\hline
62&1&77&2&92&3&107&6&122&5&137&7\\
\hline
63&2&78&1&93&2&108&3&124&5&138&5\\
\hline
64&2&79&4&94&4&109&6&124&3&139&11\\
\hline
65&0&80&2&95&3&110&0&125&6&140&5\\
\hline
66&2&81&2&96&4&111&6&126&8&141&6\\
\hline
67&1&82&2&97&2&112&6&127&6&142&8\\
\hline
68&2&83&3&98&3&113&5&128&6&143&3\\
\hline
69&2&84&5&99&2&114&3&129&7&144&8\\
\hline
70&1&85&4&100&4&115&7&130&4&145&8\\
\hline
71&2&86&4&101&5&116&6&131&9&146&7\\
\hline
72&2&87&3&102&5&117&2&132&2&147&11\\
\hline
73&1&88&2&103&5&118&6&133&8&148&5\\
\hline
74&3&89&3&104&4&119&9&134&9&149&10\\
\hline
75&3&90&2&105&4&120&8&135&5&150&6\\
\hline
\end{tabular}
\end{center}

\appendix
\section{Appendix: $d=46,\ 50,\ 52,\ 54,\ 57$, by V. Gritsenko}

In this appendix we find a vector $l\in E_8$ of square $2d$
orthogonal to exactly $12$ roots in $E_8$, where
$d$ is as in the title of the appendix.
(See \cite{K3} and \cite{sympl} for the general context of this question.)
We use below the combinatorics of the Dynkin diagram of $E_8$.
We take the Coxeter basis of simple roots in $E_8$ as in \cite{Bou}:
\noindent

\begin{picture}(300,10)(55,10)
\put(100,0){\circle*{5}}
\put(95,10){$\alpha_1$}
\put(100,0){\vector(1,0){42}}
\put(140,0){\circle*{5}}
\put(135,10){$\alpha_3$}
\put(140,0){\vector(1,0){42}}
\put(180,0){\circle*{5}}
\put(175,10){$\alpha_4$}
\put(180,1){\vector(0,-1){43}}
\put(180,-40){\circle*{5}}
\put(175,-50){$\alpha_2$}
\put(180,0){\vector(1,0){42}}
\put(220,0){\circle*{5}}
\put(215,10){$\alpha_5$}
\put(220,0){\vector(1,0){42}}
\put(260,0){\circle*{5}}
\put(255,10){$\alpha_6$}
\put(260,0){\vector(1,0){42}}
\put(300,0){\circle*{5}}
\put(295,10){$\alpha_7$}
\put(300,0){\vector(1,0){42}}
\put(340,0){\circle*{5}}
\put(335,10){$\alpha_8$}
\end{picture}
\vskip2.2cm
\noindent
where $(e_1,\dots, e_8)$ is a Euclidean basis in the lattice $\ZZ^8$ and
\begin{gather*}
\alpha_1=\frac 1{2}(e_1+e_8)-\frac 1{2}(e_2+e_3+e_4+e_5+e_6+e_7),\\
\alpha_2=e_1+e_2,\quad \alpha_k=e_{k-1}-e_{k-2}\ \ (3\le k\le 8).
\end{gather*}
The lattice $E_8$ contains $240$ roots.
We recall that any root is a sum of simple roots with integral coefficients
of the same sign.
The fundamental weights $\omega_j$ of $E_8$ form the dual basis in
$E_8=E_8^\vee$, so
$(\alpha_i,\,\omega_j)=\delta_{ij}$. The formulae for the weights
are given in \cite[Tabl. VII]{Bou}. The Cartan matrix of the dual basis is
\begin{equation}\label{Cartan-w}
((\omega_i,\,\omega_j))=\left(
\begin{matrix}
4&5&7&10&8&6&4&2\\
5&8&10&15&12&9&6&3\\
7&10&14&20&16&12&8&4\\
10&15&20&30&24&18&12&6\\
8&12&16&24&20&15&10&5\\
6&9&12&18&15&12&8&4\\
4&6&8&12&10&8&6&3\\
2&3&4&6&5&4&3&2
\end{matrix}
\right).
\end{equation}
We consider the two following cases when the orthogonal complement of
a vector $l$ in $E_8$ contains exactly $12$ roots: $
R(l_{E_8}^\perp)=A_2\oplus 3A_1$ or $A_2\oplus A_2$.  (We note that
$\# R(A_1)=2$ and $\# R(A_2)=6$.)
\medskip

\noindent{\bf I: $\bf d=46,\ 50,\ 54,\ 57$.}  There are four possible
choices of the subsystem $A_2\oplus 3A_1$ inside the Dynkin diagram of
$E_8$ according to the choices of simple roots of $A_2$, namely
$A_2^{(1,3)}=\latt{\alpha_1, \alpha_3}$, $A_2^{(2,4)}=\latt{\alpha_2,
  \alpha_4}$, $A_2^{(5,6)}=\latt{\alpha_5, \alpha_6}$ or
$A_2^{(7,8)}=\latt{\alpha_7, \alpha_8}$.  If $A_2$ is fixed then the
three pairwise orthogonal copies of $A_1$ in the Dynkin diagram are
defined automatically.

First, we consider $A_2^{(5,6)}=\latt{\alpha_5, \alpha_6}$.  Then
$3A_1^{(5,6)}=\latt{\alpha_2}\oplus \latt{\alpha_3}\oplus
\latt{\alpha_8}$.  Moreover $A_2^{(5,6)}\oplus 3A_1^{(5,6)}$ is the
root system of the orthogonal complement of the vector
$l_{5,6}=\omega_1+\omega_4+\omega_7\in E_8$.  In fact, if
$r=\sum_{i=1}^8 x_i \alpha_i$ is a positive root ($x_i\ge 0$) then
$(r,\,l_{5,6})=x_1+x_4+x_7=0$. Therefore $x_1=x_4=x_7=0$ and $r$
belongs to $A_2^{(5,6)}\oplus 3A_1^{(5,6)}$.  Using the Cartan matrix
(\ref{Cartan-w}) we obtain that $l_{5,6}^2=2\cdot 46$.  Doing similar
calculations with the other three copies of $A_2$ given above we find
$$
l_{1,3}=\omega_4+\omega_6+\omega_8,\quad
l_{2,4}=\omega_3+\omega_5+\omega_7,\quad
l_{7,8}=\omega_1+\omega_4+\omega_6
$$
with
$l_{1,3}^2=2\cdot 50$, $l_{2,4}^2=2\cdot 54$ and $l_{7,8}^2=2\cdot 57$.
\medskip

\noindent{\bf II: ${\bf d=52}$.}  We consider the sublattice
$M=A_2\oplus A_2=\latt{\alpha_3, \alpha_4}\oplus \latt{\alpha_6,
  \alpha_7}$ in $E_8$.  Then $M$ is the root system of the orthogonal
complement of the vector $l_M=\omega_1+\omega_2+\omega_5+\omega_8$
with $l_M^2=2\cdot 52$.

\medskip
\noindent
V.A.~Gritsenko, Universit\'e Lille 1, Laboratoire Paul Painlev\'e,
F-59655 Villeneuve d'Ascq, Cedex, France

\noindent
{\tt valery.gritsenko@math.univ-lille1.fr}

\section{Appendix: Computer code}

Below is the code used to check the combinatorial problem
from~\cite{K3}, and thus to find Theorem~\ref{thm52}. The programs
were written in the functional programming language Haskell ({\tt
  http://www.haskell.org}). The web page 
{\tt http://people.bath.ac.uk/masgks/Rootcounts}
contains links to further code and output.

\begin{verbatim}
{-# LANGUAGE TypeSynonymInstances,NoImplicitPrelude #-}
module E8 where

import qualified Algebra.Ring
import           Control.Applicative      ((<$>),(<*>))
import qualified Data.Vector          as V
import           Data.List                (intercalate,nubBy)
import qualified Data.MemoCombinators as Memo
import           Data.Ratio               
  (Ratio,numerator,denominator,(%))
import qualified Data.Set             as Set
import           Data.Typeable            (Typeable)
import           Math.Combinatorics.Species 
  (ksubsets,set,ofSize,enumerate,Set(getSet,Set),Prod(Prod))
import           MyPrelude hiding (numerator,denominator,(%))
import qualified Prelude
import           System.Environment       (getArgs)
import qualified Algebra.Additive

 -- Some types and helper functions for dealing with 
 -- "vectors" (implemented as arrays of rational numbers).

type Coordinate
  = Ratio Int

type Vector
  = V.Vector Coordinate

-- Inner product.
inp :: Vector -> Vector -> Coordinate
inp a b = V.sum (V.zipWith (*) a b)

half :: Coordinate
half = 1 % 2

-- Product of scalar with vector.
l :: Coordinate -> Vector -> Vector
l = V.map . (*)

instance Algebra.Additive.C Vector where
  (+) = V.zipWith (+)
  (-) = V.zipWith (-)
  negate = l (-1)
  zero = V.fromList [0,0,0,0,0,0,0,0]

-- Some data regarding E_8

delta :: (Eq a,Algebra.Ring.C b) => a -> a -> b
delta i j = if i == j then 1 else 0

-- 'e i' gives the i'th standard basis vector of R_8.
e :: Int -> Vector
e i = V.fromList $ map (delta i) [1 .. 8]

-- This is the usual integral basis of the lattice E_8.
basis :: [Vector]
basis =
  [
    l half $ (e 1 + e 8) - (sum $ map e [2 .. 7])
  , e 1 + e 2
  ] ++ map (\ i -> e (i - 1) - e (i - 2)) [3 .. 8]

roots :: [Vector]
roots = d8 ++ x118 where
  d8 = concatMap ((\ [a,b] -> 
    [a + b,a - b,b - a,negate a - b]) . map e . getSet) $
      enumerate (ksubsets 2) [1 .. 8]
  x118 = map (\ (Prod (Set neg) (Set pos)) -> 
    l half $ sum (map (negate . e) neg) + sum (map e pos)) $
      enumerate ((set `ofSize` even) * set) [1 .. 8]

-- 'posRoots' contains exactly one of every pair 
-- (a,-a) of roots.
posRoots :: [Vector]
posRoots = nubBy (\ a b -> a == b || a == negate b) roots

-- Generate elements l of the E_8 lattice with the property 
-- that l^2 = 2 d.  We need only one element of each orbit 
-- under the action of the Weyl group.  In particular, we 
-- may assume that all coordinates but one (say, the first) 
-- are nonnegative, and that the successive coordinates are 
-- nondecreasing.  We generate exactly one element of each 
-- H-orbit, where H is the subgroup of permutations and even 
-- sign changes.

gen :: Int -> [Vector]
gen d = genInt d ++ genHalfInt d

genInt :: Int -> [Vector]
genInt d = map (V.fromList . map fromIntegral) $ go [] 0 where
  -- Given the length of a partial vector, compute the maximal 
  -- new coordinate which does not increase the length of the 
  -- vector beyond 2 d.
  maxCoord :: Int -> Int
  maxCoord s = floor (sqrt (fromIntegral $ dD - s) :: Double)
  
  dD :: Int
  dD = 2 * d
  
  -- We maintain a list of coordinates chosen so far, every 
  -- one together with the sum of squares of the coordinates 
  -- up to and including that coordinate. 
  -- The generated vectors are elements of E_8, because the 
  -- sum of the squares of their components is even, hence 
  -- the sum of the components as well.
  go :: [(Int,Int)] -> Int -> [[Int]]
  -- We have fixed all eight coordinates.
  go fixed@((_,sq) : ps) 8
    -- The vector has the right length; add the relevant 
    -- solutions (using 'vary'), and continue searching.
    | sq == dD  = vary (map fst fixed) ++ lower ps 7
    -- The vector has the wrong length, continue searching.
    | otherwise = lower ps 7
  go fixed               n = let
    (m,s) = case fixed of
      []        -> (maxCoord 0,0)
      (c,s) : _ -> (Prelude.min (maxCoord s) c,s)
   in
    go ((m,s + m ^ 2) : fixed) (n + 1)
  
  -- Lexicographically decrease the given vector, and continue
  -- the generation from there.
  lower :: [(Int,Int)] -> Int -> [[Int]]
  lower []           _ = []
  lower ((x,s) : ps) n
    | x == 0    = lower ps (n - 1)
    | otherwise = go ((x - 1,s + 1 - 2 * x) : ps) n
  
  vary :: [Int] -> [[Int]]
  vary (x : xs) = if x == 0
    then [0 : xs]
    else [x : xs,negate x : xs]

-- For vectors with all coordinates half-integers, we work 
-- with the doubles of the coordinates.
genHalfInt :: Int -> [Vector]
genHalfInt d = map (V.fromList . map (% 2)) $ go [] 0 where
  maxCoord :: Int -> Int
  maxCoord = Memo.integral m where
    m s = f $ floor (sqrt (fromIntegral $ dE - s) :: Double)
    f k = if odd k then k else k - 1
  
  dE :: Int
  dE = 8 * d
  
  go :: [(Int,Int)] -> Int -> [[Int]]
  go fixed@((_,sq) : ps) 8
    | sq == dE  = filter e8 (vary $ map fst fixed) 
               ++ lower ps 7
    | otherwise = lower ps 7
  go fixed                   n = let
    (m,s) = case fixed of
      []        -> (maxCoord 0,0)
      (c,s) : _ -> (Prelude.min (maxCoord s) c,s)
   in
    go ((m,s + m ^ 2) : fixed) (n + 1)
  
  -- Decides whether a given vector is an element of E_8
  e8 :: [Int] -> Bool
  e8 = (== 0) . flip rem 4 . sum
  
  lower :: [(Int,Int)] -> Int -> [[Int]]
  lower []           _ = []
  lower ((x,s) : ps) n
    | x == 1    = lower ps (n - 1)
    | otherwise = go ((x - 2,s + 4 - 4 * x) : ps) n

  vary :: [Int] -> [[Int]]
  vary (x : xs) = [x : xs,negate x : xs]
\end{verbatim}

\bibliographystyle{alpha}

\bigskip

\noindent
A. Peterson, Korteweg de Vries Instituut voor Wiskunde,
Universiteit van Amsterdam, P.O. Box 9424, 1090 GE Amsterdam, The
Netherlands

\noindent
{\tt ariep@xs4all.nl}

\medskip
\noindent
G.K.~Sankaran, Department of Mathematical Sciences, University of
Bath, Bath BA2 7AY, England

\noindent
{\tt gks@maths.bath.ac.uk}

\end{document}